\documentclass{amsart}

\usepackage{amssymb,amsfonts, latexsym,amsmath,hhline,array,longtable}
\usepackage{pdfsync,color, comment,colortbl, mathrsfs,stmaryrd,cite,graphicx,yfonts,enumerate}

\usepackage{ifpdf}
\ifpdf \usepackage[colorlinks=true, citecolor=blue, linkcolor=blue, urlcolor=blue]{hyperref} \fi

\newcommand{\cal}{\mathcal}

\newtheorem{formula}{}[section]
\newtheorem{definition}[formula]{Definition}
\newtheorem{corollary}[formula]{Corollary}
\newtheorem{remark}[formula]{Remark}
\newtheorem{lemma}[formula]{Lemma}
\newtheorem{theorem}[formula]{Theorem}
\newtheorem*{claim}{Claim}

\def\thrm{\begin{theorem}}
\def\thrml#1{\begin{theorem}\label{#1}}
\def\ethrm{\end{theorem}}
\def\rmrk{\begin{remark}}
\def\rmrkl#1{\begin{remark}\label{#1}}
\def\ermrk{\end{remark}}
\def\dfntn{\begin{definition}}
\def\dfntnl#1{\begin{definition}\label{#1}}
\def\edfntn{\end{definition}}
\def\nmrt{\begin{enumerate}}
\def\enmrt{\end{enumerate}}
\def\tm#1{\item[{\rm (#1)}]}
\def\qtnl#1{\begin{equation}\label{#1}}
\def\eqtn{\end{equation}}
\def\lmm{\begin{lemma}}
\def\lmml#1{\begin{lemma}\label{#1}}
\def\elmm{\end{lemma}}
\def\crllr{\begin{corollary}}
\def\crllrl#1{\begin{corollary}\label{#1}}
\def\ecrllr{\end{corollary}}
\def\css{\begin{cases}}
\def\ecss{\end{cases}}
\def\prf{\begin{proof}}
\def\eprf{\end{proof}}
\def\clm{\begin{claim}}
\def\eclm{\end{claim}}

\def\cA{{\cal A}}

\def\cG{{\cal G}}

\def\cR{{\cal R}}

\def\mC{{\mathbb C}}
\def\mF{{\mathbb F}}

\def\fX{{\mathfrak X}}

\DeclareMathOperator{\aut}{Aut}

\DeclareMathOperator{\AGL}{AGL}
\DeclareMathOperator{\AGaL}{A{\rm \Gamma}L}
\DeclareMathOperator{\cay}{Cay}

\DeclareMathOperator{\GaL}{{\rm \Gamma}L}
\DeclareMathOperator{\GamO}{\Gamma O}

\DeclareMathOperator{\GL}{GL}

\DeclareMathOperator{\id}{id}

\DeclareMathOperator{\ord}{ord}

\DeclareMathOperator{\SL}{SL}

\DeclareMathOperator{\sym}{Sym}
\DeclareMathOperator{\VO}{VO}
\DeclareMathOperator{\VSz}{VSz}
\DeclareMathOperator{\VD}{VD}

\def\qaq{\quad\text{and}\quad}

\def\mltst#1{[#1]}
\def\ov{\overline}

\begin{document}

%\begin{comment}
\title{On multivalued groups of order~$3$}
\author{Jin Guo}
\address{School of Mathematics and Statistics, Hainan University, Haikou {\rm570228}, Hainan, P.~R.~China}
\email{guojinecho@163.com}
\author{Ilia Ponomarenko}
\address{Steklov Institute of Mathematics at St. Petersburg, Russia}
\email{inp@pdmi.ras.ru}
\author{Andrey V. Vasil'ev}
\address{Sobolev Institute of Mathematics, Novosibirsk 630090, Russia}
\email{vasand@math.nsc.ru}
%\thanks{J. Guo and I. Ponomarenko were supported by the National Natural Science Foundation of China (Grant No. 12361003), A.~V.~Vasil'ev was supported by the Sobolev Institute of Mathematics state contract (FWNF-2022-0002) and by the National Natural Science Foundation of China (Grant No. 12171126).}
\date{}

\begin{abstract}
A complete classification of the multivalued coset groups of order~$3$ is given. The proof is based on the classification of rank~$3$ groups having regular normal subgroups.


\end{abstract}

\maketitle

\section{Introduction}
Given a positive integer $n$, an \emph{$n$-valued group} on a set $X$  can be thought as a natural generalization of an ordinary group,  where the product of any two elements  is a multiset containing   elements of~$X$ with some multiplicities  the sum of which  is equal to~$n$ (as for the exact definitions, we refer to Section~\ref{170824d}). Thus, a $1$-valued group is nothing else than a group in the usual sense. The concept of multivalued group goes back to a construction in the Buchshtaber-Novikov paper~\cite{BuchN1971}. Nowadays, the theory of multivalued groups is an actively developing part of mathematics and finds numerous applications in its various fields, such as topology~\cite{BukH1982,BukР1991,Bu2012}, dynamical systems~\cite{BuV1996,BuV2019,GaiYag07}, Hopf algebras~\cite{BuR1997,BuR1998},
discrete mathematics~\cite{BuVEP1996e,Yag2005} and so on, see also the survey \cite{Bu2006} and references therein. Let us specially note the recent completion of the classification of involutive $2$-valued groups in~\cite{BuVG2022,Gaif2024}.

In the general definition of a multivalued group, the operation of taking the inverse is not necessarily an involutive antiisomorphism of its group algebra (as in the case of an ordinary group). If a property of this type is postulated from the beginning, then the group algebras of the corresponding multivalued groups appear in a natural correspondence with so called combinatorial algebras introduced in~\cite{BuVEP1996e}. The finite dimensional combinatorial  algebras form a natural subclass of generalized table algebras used to study the characters of finite groups~\cite{AradFM1999}. The multivalued groups whose  group algebras are combinatorial are said to be \emph{involutive}.\footnote{Here we use the definition of involutive group from \cite{BuVEP1996e} which is weaker than that in~\cite{BuVG2022}.} A direct definition of multivalued involutive  groups (not using group algebras) is given in Section~\ref{170824d}.

A natural example of multivalued involutive groups is represented by the \emph{coset groups}, see~\cite{BuR1997}. Let $G$ and $A$ be ordinary ($1$-valued) groups, and let $\varphi:G\to A$ be a homomorphism. Then the triple $(G,A,\varphi)$ yields an $|A|$-valued group on the set~$X$ of the orbits of the action of $A$ on~$G$, afforded by~$\varphi$. %In particular, the coset group  is finite whenever the number $|X|$ is finite.

In Section~\ref{100824b}, we present a construction of multivalued involutive groups of order~$3$ from \emph{strongly regular graphs}; concerning a huge and nice theory of these graphs, we refer the reader to recent monograph~\cite{BrouwerM2022}, see also Section~\ref{100824b}. Here we only note that each strongly regular graph $\cG$ corresponds to a triple of integer parameters $k>0$ and $\lambda,\mu\ge 0$. The multivalued involutive group corresponding to such~$\cG$  is denoted by $\fX(\cG)=\fX(k,\lambda,\mu)$. Among the strongly regular graphs, the \emph{rank $3$ graphs} are of   special interest (in fact, the most part of the monograph~\cite{BrouwerM2022} is devoted to complete classification of them). The defining property of such a graph $\cG$ is that the group $\aut(\cG)$ is of rank~$3$, i.e., is  transitive on the vertices of~$\cG$,  on the pairs of adjacent vertices of~$\cG$, and on the  pairs of adjacent vertices of the complement of~$\cG$.

%Let us proceed to the formulation of the main result of the present paper.
It is easily seen that there is a unique (up to isomorphism)  multivalued group of order~$2$. Our goal  is to  classify the coset multivalued groups of order~$3$. One family of these groups  is given by the $(2k+1)$-valued groups $X(k)$, where $4k+3$ is a prime power (see~\cite{Ponomarenko2024a} and the last paragraph of Section~\ref{170824d}). One more  family consists of the multivalued groups $\fX(\cG)$, where $\cG$ is a rank~$3$ graph. It turns out that these two families exhaust the whole variety of coset multivalued groups of order~$3$.

\thrml{100824a}
A multivalued  group $\fX$ of order~$3$ is a coset group if and only if $\fX$ is isomorphic to a multivalued group $X(k)$, where $4k+3$ is a prime power, or a multivalued group $\fX(k,\lambda,\mu)$, where the integers $k,\lambda,\mu$ are as in  Theorem~{\rm\ref{180824c}}.
\ethrm

The idea of the proof of Theorem~\ref{100824a} is based on the fact that if $\fX$ is a multivalued  group associated with triple $(G,A,\varphi)$, then the semidirect product $G\rtimes \varphi(A)\le \sym(G)$ is a group of rank~$3$. Since  $G$ is a normal regular subgroup of the semidirect product, the problem reduces to determine the rank $3$-groups having a normal regular subgroup. At this point, we prove that this subgroup must be elementary abelian. This enables us to  use the classification of the affine rank~$3$-groups, completed in~\cite{liebeckAffine}, see also  \cite[Chapter~11]{BrouwerM2022}.

All undefined terms concerning multivalued groups,  permutation groups, and strongly regular graphs can be found in~\cite{Bu2006}, \cite{DixM1996}, and \cite{BrouwerM2022}, respectively.

The authors thank V.~M. Buchshtaber   for fruitful discussions and useful comments.

\section{Multivalued groups}\label{170824d}

Let $X$ be a set,  and let $n$ be a positive integer. A \emph{$n$-valued group} $\fX= (X, \cdot, \star)$ on~$X$ is defined by the operations $\cdot$  and  $\star$ of multiplication and  taking  inverse, respectively.  The product of  two elements $x,y\in X$ is defined to be  an $n$-multiset
$$
 x\cdot y=\mltst{(x\cdot y)_1,\ldots,(x\cdot y)_n}
$$
with elements belonging to~$X$. It is assumed that the multiplication is associative, i.e., the $n^2$-multisets $(x\cdot y)\cdot z$ and $x\cdot(y\cdot z)$ coincide for all $x,y,z\in X$,
\qtnl{160824a}
\mltst{((x\cdot y)_i\cdot z)_j:\ 1\le i,j\le n}=\mltst{(x\cdot (y\cdot z)_i)_j:\ 1\le i,j\le n}.
\eqtn
Furthermore, it is assumed that there is an \emph{identity} element $e\in X$ with respect to the multiplication, i.e., for all $x\in X$,
$$
(x\cdot e)_i=(e\cdot x)_i=x,\qquad i=1,\ldots,n.
$$
Finally, the  operation of taking inverse maps $X$ to $X$ so that $e$ is an element of the multisets $x\cdot x^\star$ and $x^\star\cdot x$. The number $|X|$ is called the \emph{order} of~$\fX$.

In what follows,  for any $x,y,z\in X$, we denote by~$m_{x,y}^z$ the multiplicity of the element $z\in X$ in the multiset $x\cdot y$.  A bijection $f: X\to X'$, $x\mapsto x'$, is called an \emph{isomorphism} of the $n$-valued group $\fX$ on $X$ onto $n'$-valued group $\fX'=(X',\cdot,\star)$  if
$$
{m_{x^{},y^{}}^{z^{}}\over{n^{}}}={m_{x',y'}^{z'}\over{n'}}\quad\text{ for all } x,y,z\in X.
$$

A \emph{coset group} is given by a pair of ordinary groups $G$ and  $A$, and a homomorphism $\varphi:A\to\aut(G)$. In the corresponding $n$-valued group $\fX(G,A,\varphi)$, the number $n$ is equal to $|A|$, the set $X$ consists of the $A$-orbits in the induced action of~$A$  on~$G$, the multiplicities are defined by
\qtnl{180824l}
m_{x,y}^z=|\{a\in A:\ g\cdot h^a\in z\}|,\qquad g\in x,\ h\in y,
\eqtn
where $h^a:=h^{\varphi(a)}$, and the inverse of $x$ is equal to the $A$-orbit $(g^{-1})^A$. In particular, any group $G$ can be treated as the $1$-valued coset group $\fX(G,1,\id)$. Every coset group is \emph{involutive}, i.e., the mapping $X\to X$, $x\mapsto x^\star$, is an involution respecting the diagonal multiplicities $m(x)=m_{x,x^\star}^e$, $x\in X$. More exactly, for all $x,y\in X$, we have
\qtnl{180824v1}
m_{x,y}^e>0\qquad \Leftrightarrow\qquad y=x^\star
\eqtn
and
\qtnl{180824v2}
m(x)=m(x^\star), \qquad  (x\cdot y)^\star = y^\star \cdot x^\star,
\eqtn
where  $(x\cdot y)^\star=\mltst{(x\cdot y)_1^\star,\ldots,(x\cdot y)_n^\star}$.

The associativity implies that  the following identity holds in any multivalued  involutive group:
\qtnl{180824q}
%m(z)\cdot m_{x,y}^{z^\star}=
m(x)\cdot m_{y,z}^{x^\star}=m(y)\cdot m_{z,x}^{y^\star},\qquad x,y,z\in X.
\eqtn
To prove the identity, we  calculate the multiplicities of the identity element $e$ in the multisets on the left- and right-hand side of formula~\eqref{160824a}.  From formula~\eqref{180824v1}, it follows that they are equal to
$$
m_{y,z}^{x^\star}\cdot m_{x^\star,x}^e=m_{y,z}^{x^\star}\cdot m(x^\star)\qaq
m_{z,x}^{y^\star}\cdot m_{y^\star,y}^e=m_{z,x}^{y^\star}\cdot m(y^\star),
$$
respectively. This proves the required identity by the first equality in formula~\eqref{180824v2}.

In the present paper, we will focus on multivalued involutive groups of order~$3$. In this case, only three multiplicities $m_{x,y}^z$ are enough to determine the group up to isomorphism.

\thrml{170824a}
Let $\fX=(X,\cdot,\star)$ be  an $n$-valued  involutive group of order~$3$, and let $X=\{e,x,y\}$, $m(x)=m_{x,x^*}^e$,  $m(y)=m_{y,y^*}^e$, and $a(x)=m_{x,x}^x$. Then exactly one of the two following statements holds:
\nmrt
\tm{i} $x^\star=x$, $y^\star=y$, and
$$
x\cdot x  = \mltst{\underbrace{e,\ldots,e}_{m(x)},\underbrace{x,\ldots,x}_{a(x)},\underbrace{y,\ldots,y}_{n-m(x)-a(x)}},\qquad
y\cdot y =\mltst{\underbrace{e,\ldots,e}_{m(y)},\underbrace{x,\ldots,x}_{a(y)},\underbrace{y,\ldots,y}_{n-m(y)-a(y)}},
$$
$$
x\cdot y=y\cdot x= \mltst{\underbrace{x,\ldots,x}_{a(x,y)},\underbrace{y,\ldots,y}_{n-a(x,y)}},
$$
where $a(x,y)=r(n-m(x)-a(x))$  and $a(y)=ra(x,y)$ with $r=m(y)/m(x)$,
\tm{ii} $x^\star=y$, $y^\star=x$, and
$$
x\cdot x = \mltst{\underbrace{x,\ldots,x}_{a(x)},\underbrace{y,\ldots,y}_{n-a(x)}},\qquad
y\cdot y =\mltst{\underbrace{x,\ldots,x}_{n-a(x)},\underbrace{y,\ldots,y}_{a(x)}},
$$
$$
x\cdot y=y\cdot x= \mltst{\underbrace{e,\ldots,e}_{n-2a(x)},\underbrace{x,\ldots,x}_{a(x)},\underbrace{y,\ldots,y}_{a(x)}}.
$$
\enmrt
\ethrm
\prf
We consider two cases depending on whether the involution $\star$ %\in\sym(3)$
is trivial or not. First, let $x^\star=x$, $y^\star=y$.  By formula~\eqref{180824v1}, we have $m_{y,z}^e=m_{z,y}^e=0$ for all $z\ne y$. Furthermore, $x\cdot y=x^\star\cdot y^\star=(y\cdot x)^\star=y\cdot x$ by the second equality in formula~\eqref{180824v2}. Thus it suffices to verify that $m_{y,y}^x=a(y)$ and $m_{x,y}^x=a(x,y)$. By the identity~\eqref{180824q}, we have
$$
m_{x,y}^x=\frac{m(y)}{m(x^\star)}\cdot m_{x^\star,x}^{y^\star}=
\frac{m(y)}{m(x)}\cdot m_{x,x}^y=
r\cdot(n-m(x)-a(x))=a(x,y),
$$
where we used the fact that $m(x)>0$ by the equivalence~\eqref{180824v1}. Using the obtained equality, we similarly  conclude that
$$
m_{y,y}^x=\frac{m(y)}{m(x^\star)}\cdot m_{x^\star,y}^{y^\star}=\frac{m(y)}{m(x)}\cdot m_{x,y}^y=r\cdot a(x,y)=a(y),
$$
as required.

Let $x^\star=y$, $y^\star=x$.  By formula~\eqref{180824v1}, we have $m_{x,x}^e=m_{y,y}^e=0$, whereas by the first and second equalities  in formula~\eqref{180824v2}, we have $m(x)=m(x^\star)=m(y)$ and $m_{y,y}^y=m_{x,x}^x=a(x)$, respectively. Thus it suffices to verify that $m_{x,y}^x=m_{x,y}^y=n-a(x)$. Using the identity~\eqref{180824q}, we have
$$
m_{x,y}^x=\frac{m(x)}{m(x^\star)}\cdot m_{y,x^\star}^{x^\star}=
m_{y,y}^y=a(x),
$$
as required.
\eprf

The $n$-valued involutive  groups of order~$3$, defined in statements (i) and (ii) of Theorem~\ref{170824a}, will be denoted by $\fX_n(m_1,m_2,a)$ with $m_1=m(x)$,  $m_2=m(y)$, $a=a(x)$,  and  $\fX_n(a)$, respectively.

\crllrl{190824j}
The multivalued groups $\fX_n(m_1,m_2,a)$ and $\fX_{n'}(m'_1,m'_2,a')$ {\rm (}respectively, $\fX_n(a)$ and $\fX_{n'}(a')${\rm)} are isomorphic if and only if $u/n=u'/n'$ for $u\in\{m_1,m_2,a\}$ {\rm(}respectively, for $u=a${\rm)}.
\ecrllr

As it is readily seen from the definition of the multivalued group $X(k)$ in~\cite{Ponomarenko2024a}, in the above notation, $X(k)=\fX_{2k+1}(k)$.

\section{Multivalued groups  and strongly regular graphs}\label{100824b}
The main goal of this section is to define a large family of multivalued groups of order~$3$, related with strongly regular graphs. %~\cite{BrouwerM2022}.
To this end,  let $\cG$ be a \emph{strongly regular} graph;  this means that $\cG$ is not complete or edgeless, and the number of common neighbors of two arbitrary vertices $\alpha$ and $\beta$ in $\cG$ is equal to $k$ (respectively, $\lambda$, $\mu$) depending on whether $\alpha$ and $\beta$ are equal (respectively, adjacent, distinct and nonadjacent). The integers $k>0$ and $\lambda, \mu\ge 0$ are called the \emph{parameters} of~$\cG$. The parameters satisfy the relation
$$
k(k-1-\lambda)=(v-k-1)\mu,
$$
where $v$ is the number of vertices of $\cG$. Because of this relation, we do not include~$v$  as a parameter of~$\cG$, as is usually done.

Let $A$ be the adjacency matrix of the graph $\cG$, that is $A$ is a \{0,1\}-matrix of size $v\times v$, where the entry $A_{\alpha,\beta}$ is equal to $1$ or $0$ depending on whether the vertices $\alpha$ and $\beta$ are adjacent or nonadjacent in~$\cG$. The algebra $\cA=\cA(\cG)$ generated over~$\mC$ by~$A$, has a linear basis consisting of the identity matrix $A_0=I$, the matrix $A_1=A$, and the adjacency matrix $A_2$ of the  complement  graph $\bar\cG$ of $\cG$. Moreover,
\begin{align}
A_1\cdot A_1 & = k\cdot A_0+\lambda\cdot A_1+\mu\cdot A_2\\ \label{190824a}
A_2\cdot A_2 &= \bar k\cdot A_0+\bar\mu\cdot A_1+\bar\lambda \cdot A_2\\ \label{190824b}
A_1\cdot A_2 =A_2\cdot A_1 & = {{v-2-\lambda-\bar\mu}\over{2}}\cdot A_1 + {{v-2-\bar\lambda-\mu}\over{2}}\cdot A_2,
\end{align}
where $\bar k= v-k-1$, $\bar\lambda =v-2k+\mu-2$, and $\bar\mu=v-2k-\lambda$. In fact, the above equalities imply that the complement graph $\bar\cG$ is strongly regular with parameters~$(\bar k, \bar\lambda,\bar\mu)$.

The algebra $\cA$ is  a \emph{combinatorial algebra} in the sense of~\cite[Definition~1.6]{BuVEP1996e} with respect to the linear basis  $\cR=\{A_0,A_1,A_2\}$ and the identical antiisomorphism~$*$.  Indeed, let $c_{r,s}^t$ be the coefficient at the matrix $A_t$ in the product $A_r\cdot A_s$, where $r,s,t\in\{0,1,2\}$. The nonnegative integers $c_{r,s}^t$ are  the structure constants of the algebra~$\cA$ with respect to the basis $\cR$. The diagonal numbers $d_r=d(A_r)=c_{r,r}^0$ are obviously bounded by~$v$ from above. Moreover,
$$
c_{r,s}^0=d_r\,\delta_{r,s},
$$
where $\delta$ is the Kronecker delta and  the mapping $A_r\mapsto d_r$, $r=0,1,2$, is a one-dimensional representation of $\cA$.%Thus, as a consequence of  Theorem~2.2 of that paper, we arrive at the following statement.

\thrml{110824a}
Let $\cG$ be a strongly regular graph with parameters $(k,\lambda,\mu)$, and let~$n$ be the least common multiple of $k$ and~$\bar k$. Then there is an $n$-valued  involutive group $\fX(\cG)=\fX(k,\lambda,\mu)$ on the set $\{x_0, x_1, x_2\}$, such that
\qtnl{120824a}
m_{x_r,x_s}^{x_t}=n\cdot c_{r,s}^t\cdot{d_t\over{d_rd_s}}, \qquad r,s,t\in\{0,1,2\}.
\eqtn
In particular,
\qtnl{1908824c}
m(x_1)=m_{x_1,x_1}^0=\frac{n}{k},\qquad
m(x_2)=m_{x_2,x_2}^0=\frac{n}{\bar k},\qquad a(x_1)=m_{x_1,x_1}^{x_1}=\frac{n\lambda}{k}.
\eqtn
\ethrm
\prf
By \cite[Theorem~2.2]{BuVEP1996e}, the combinatorial algebra $\cA=\cA(k,\lambda,\mu)$ defines an $n'$-valued involutive  group $\fX'$ on the set  $\{x_0, x_1, x_2\}$ with  multiplication  defined by formula~\eqref {120824a} in which the number $n$ is replaced by the least common multiple~$n'$ of the (positive) denominators in the representation of the rational number $c_{r,s}^t/d_r$ in the form of a reduced fraction.\footnote{There is a misprint in formula~(7) of the paper~\cite{BuVEP1996e}: the great common divisor there should be replaced by the least common multiple.} Since $n'$ divides $n$, the multivalued group $\fX'$ is isomorphic to an $n$-valued  involutive group $\fX(k,\lambda,\mu)$ with multiplication  defined by formula~\eqref {120824a}.
\eprf

{\bf Example.} The Petersen graph is a unique strongly regular graph on $10$ vertices  with parameters $(3,0,1)$. Thus,   $\fX(3,0,1)$ is a $6$-valued group on the $3$-element set $X=\{x_0,x_1,x_2\}$ and the multiplication  defined as follows:
$$
x_1\cdot x_1 = [\underbrace{x_0,\ldots,x_0}_{2},\underbrace{x_2,\ldots,x_2}_{4}],\qquad
x_2\cdot x_2 = [x_0,\underbrace{x_1,\ldots,x_1}_{2},\underbrace{x_2,\ldots,x_2}_{3}],
$$
$$
x_1\cdot x_2=x_2\cdot x_1  = [\underbrace{x_1,\ldots,x_1}_{2},\underbrace{x_2,\ldots,x_2}_{4}].
$$
In the notation of Section~\ref{170824d}, we have $m_1=m(x_1)=2$, $m_2=m(x_2)=1$, and $a=a(x_1)=0$. Thus, $\fX(3,0,1)=\fX_6(2,1,0)$.\medskip

It should be noted that not each multivalued  involutive group of order $3$ is of the form $\fX(\cG)$ for some strongly regular graph~$\cG$. Indeed, %as was observed in~\cite{Ponomarenko2024a},
a straightforward computation shows that for each natural number $q = 4\ell+ 1$, there exists a multivalued group $\fX=\fX_{2\ell}(1,1,\ell-1)$. Assume that $\fX=\fX(\cG)$ for some strongly regular graph~$\cG$. Then from formulas~\eqref{120824a}, one can find that the parameters $v$ and $k$ of $\cG$ are equal to $4\ell+1$ and $2\ell$, respectively. It follows that $\bar k=2\ell=k$. As is well known, this is possible only if $v$ is a sum of two squares (see, e.g., \cite[Theorem~8.2.3]{BrouwerM2022}). Thus $\fX=\fX(\cG)$ only if $4\ell+1$ is a sum of two squares.

%\rmrkl{110824n}
The multivalued group $\fX(\cG)$ defined by a strongly regular graph $\cG$  depends on the parameters of $\cG$ only. Therefore, any two nonisomorphic strongly regular graphs with the same parameters define the same multivalued group. Note that for infinitely many positive integers $v$ there can be exponentially many  nonisomorphic strongly regular graphs with $v$ vertices and the same parameters~\cite[Subsection~1.1.18]{BrouwerM2022}.
%\ermrk

Among strongly regular graphs, the \emph{rank $3$ graphs} are of special interest. By definition these are graphs admitting an automorphism group acting transitively on the  vertices, on the ordered pairs of adjacent vertices, and on the ordered pairs of nonadjacent vertices. The automorphism group of a rank~$3$ graph is always a \emph{rank~$3$ group}, i.e., a transitive group a point stabilizer of which has exactly three orbits.

\thrml{170824q}
Let $K$ be a group of rank~$3$, $G$  a normal regular subgroup of $K$, and $A$ the stabilizer of the identity element~$e$ of~$G$ in~$K$. Assume that $A$ is of even order and $x$ is an $A$-orbit other than $\{e\}$. Then the Cayley graph $\cG=\cay(G,x)$ is  a rank~$3$ graph. Moreover, the multivalued groups $\fX(\cG)$   and  $\fX(G,A,\id)$ are isomorphic.
\ethrm
\prf
Since $A$ is of even order, we have $x^{-1}=x$. Consequently the Cayley graph~$\cG$ is undirected and hence of rank~$3$. Denote its parameters by $(k,\lambda,\mu)$. From  formulas \eqref{190824a} and~\eqref{190824b} for the intersection numbers of the combinatorial algebra $\cA(\cG)$ corresponding to the strongly regular graph~$G$, we obtain
$$
c_{1,1}^0=d_1=k,\quad  c_{2,2}^0=d_2=\bar k,\quad c_{1,1}^1=\lambda.
$$
%where $v=|G|$ is the number of vertices of $\cG$.
By Theorem~\ref{110824a}, this gives the following formulas for  the corresponding  multiplicities of the multivalued  group $\fX'=\fX(\cG)$ on the set $\{e', x', y'\}$:
\qtnl{1908824c11}
m(x')=\frac{n'}{k},\qquad
m(y')=\frac{n'}{\bar k},\qquad
m_{x',x'}^{x'}=\frac{n'\lambda}{\bar k},
\eqtn
where $n'$ is the least common multiple of $k$ and $\bar k$.

The multivalued group $\fX=\fX(G,A,\id)$ is defined on the set $\{e, x, y\}$, where we identify $e$ with the singleton $\{e\}$ and denote by $y$  the set of all nonidentity elements of the group~$G$, not belonging to~$x$. Since $x=x^{-1}$, Theorem~\ref{170824a} implies that
$$
\fX=\fX_n(m(x),m(y),a(x)),
$$
where $n=|A|$ and $a(x)=m_{x,x}^x$. For every $g,h\in x$ and $a\in A$, we have $g\cdot h^a=e$ if and only if $h^a=g^{-1}$. By formula \eqref{180824l}, this shows that  the number $m_{x,x}^e$ is equal to the order of  the stabilizer $A_h$ of the element~$h$ in the group~$A$. Since $x$ is an $A$-orbit and the cardinality of $x$ is equal to~$k$, we conclude that
\qtnl{1908824x}
m(x)=m_{x,x}^e=|A_h|=\frac{|A|}{|x|}=\frac{n}{k}.
\eqtn
Similarly, for any element $h\in y$, we have
\qtnl{1908824x1}
m(y)=m_{y,y}^e=|A_h|=\frac{|A|}{|y|}=\frac{n}{\bar k}.
\eqtn
Finally, fix  $g\in x$. Again by formula \eqref{180824l}, the number $m_{x,x}^x$ is equal to the number of all $a\in A$ for which $g\cdot (g^{-1})^a\in x$ or, equivalently, $g^a\cdot  g^{-1}\in x^{-1}=x$.  When $a$ runs over $A$, the vertex $ g^a$ runs over all the vertices of $\cG$ (each counted $|A_g|$ times) adjacent with both the vertex~$e$ (because $ g^a\in x$) and $g$  (because $ g^a\cdot  g^{-1}\in x$). Thus,
\qtnl{1908824x2}
m_{x,x}^x=|A_g|\cdot\lambda=\frac{|A|}{|x|}\lambda=\frac{n\,\lambda}{k}.
\eqtn
Now comparing formulas \eqref{1908824c} with formulas  \eqref{1908824x}, \eqref{1908824x1}, \eqref{1908824x2}, we obtain
$$
\frac{m(x)}{n}=\frac{m(x')}{n'},\quad \frac{m(y)}{n}=\frac{m(y')}{n'},\quad \frac{m_{x^{},x^{}}^{x^{}}}{n}=\frac{m_{x',x'}^{x'}}{n'}.
$$
Thus the multivalued groups $\fX$ and $\fX'$ are isomorphic by Corollary~\ref{190824j}.
\eprf

\section{The groups of rank~$3$ with regular normal subgroup}\label{280824a}

Recall that every rank $3$ group $K\leq\sym(\Omega)$ of even order corresponds to a pair of complementary rank $3$ graphs with vertex set $\Omega$. By the parameters of $K$, we mean the parameters $k$, $\lambda$, $\mu$ of the one of these graphs whose valency $k$ is lower. We denote this graph by $\mathcal{G}(K)$.

Let us emphasize two things. First, nonisomorphic rank $3$ groups $K_1,K_2\leq\sym(\Omega)$ can correspond to the same graph $\mathcal{G}(K_1)=\mathcal{G}(K_2)$. Second, nonisomorphic graphs $\mathcal{G}(K_1)$ and $\mathcal{G}(K_2)$ might have the same parameters $k$, $\lambda$, $\mu$. The goal of this section is to find all attainable triples of the parameters for the rank $3$ groups with regular normal subgroups. For every such triple $(k,\lambda,\mu)$, we also point out a possible rank $3$ graph $\mathcal{G}$ and possible rank $3$ group $K$ such that $\mathcal{G}=\cG(K)$ has the parameters $k$, $\lambda$, $\mu$.

A permutation group $K \leq \sym(\Omega)$ is called {\em affine}, if it has a regular normal elementary abelian $p$-subgroup $V$. Then $K =V:A$ is a split extension of $V$ and a point stabilizer $A$. The set $\Omega$ can be identified with $V$ considered as a vector space over the prime field $\mF_p$ in such a way that $V$ itself acts on this space by translations, and $A$ being the stabilizer of the zero vector acts on it as a subgroup of $\GL(V)$. Such an identification provides a natural embedding $K\leq\operatorname{AGL}(V)$. %\inp{Undefined notation concerning graphs and groups of rank $3$ from the theorem below can be found in~\cite{BrouwerM2022}.}

\thrml{180824c}
Let $K$ be a rank $3$ group of even order. Assume that $K$ has a regular normal subgroup $V$ of order~$v$. Then $K=V:A$ is affine, $v$ is a power of a~prime $p$, and the parameters of $K$ are either in Table~{\rm\ref{param2}} or listed below:
\nmrt[{\rm(i)}]
\item $v=p^{t+s}$, where $t,s>0$, $k=p^t-1$, $\lambda=p^t-2$, $\mu=0$,

\noindent $K$ is imprimitive, $\cG(K)$ is a disjoint union of $p^s$ cliques of size $p^t$, and as $A$ can be taken the stabilizer in $\GL(V)$ of a $t$-dimensional subspace of $V;$

\item $v=q^2$, $k=2(q-1)$, $\lambda=q-2$, $\mu=2$,

\noindent $K=K_0\wr\sym(2)$ preserves a product decomposition $V=U\times U$, where $U$ is a normal regular subgroup of $K_0$ of order $q$, $K_0$ is transitive on $U^{\#}$, $\cG(K)$ is a $q\times q$-grid, and as a point stabilizer of $K_0$ can be taken the group $\GL_1(U)${\rm;}

\item $v=4t+1$, $k=2t$, $\lambda=t-1$, $\mu=t$,

\noindent $K\le\AGaL_1(v)$, and as $\cG(K)$ and as $A$ can be taken the Paley graph and the subgroup of index $2$ in $\GL_1(V)$, respectively{\rm;}

\item  $v=p^{(e-1)t}$, where $e>2$ is a prime with  $\ord_e(p)=e-1$\footnote{to avoid coincidences with the parameters of groups from other series let $(p,e,t)\neq(2,3,2),(5,3,1),(2,3,3),(3,5,1),(2,5,2),(3,7,1),(2,11,1),(2,13,1)$.}, $k=(v-1)/e$, $\lambda=(v-3e+1-(-1)^t(e-2)(e-1)\sqrt{v})/e^2$, $\mu=(v-e+1+(-1)^t(e-2)\sqrt{v})/e^2$,

\noindent $K\le\AGaL_1(v)$, $\cG(K)$ is the Van Lint--Schrijver graph, and as $A$ can be taken a subgroup of index $e$ in $\GaL_1(V);$

\item  $v=q^{2e}$ with $e\geq3$, $k=(q+1)(q^e-1)$, $\lambda=q^e+(q-2)(q+1)$, $\mu=q(q+1)$,

\noindent $K\leq\AGaL_{2e}(q)$, $\cG(K)$ is the bilinear forms graph $H_q(2,e)$, and as $A$ can be taken the subgroup $(\GL_2(q)\circ\GL_m(q)):\aut(\mF_q);$

\item $v=q^{2e}$, $e\geq2$, $\varepsilon=\pm$, and $(q,\varepsilon)\neq(2,+)$, $k=(q^e-\varepsilon1)(q^{e-1}+\varepsilon1)$, $\lambda=q(q^{e-1}-\varepsilon1)(q^{e-2}+\varepsilon1)+q-2$, $\mu=q^{e-1}(q^{e-1}+\varepsilon1)$,

\noindent $K\leq\AGaL_{2e}(q)$ $\cG(K)$ is the affine polar graph $\VO_{2e}^{\varepsilon}(q)$, and as $A$ can be taken the subgroup $\GamO^{\varepsilon}_{2e}(q);$

\item $v=2^{2e}$, $e\geq2$, $k=2^{e-1}(2^e-1)$, $\lambda=\mu=2^{e-1}(2^{e-1}-1)$,

\noindent $K\leq\AGaL_{2e}(2)$, $\cG(K)$ is the complement $\ov{\VO_{2e}^{+}(2)}$ of the affine polar graph, and as $A$ can be taken the subgroup $\GamO^{+}_{2e}(2);$

\item $v=q^{10}$, $k=(q^2+1)(q^5-1)$, $\lambda=q^5+q^4-q^2-2$, $\mu=q^2(q^2+1)$,

\noindent $K\leq\AGaL_{10}(q)$, $\cG(K)$ is the alternating forms graph $A(5,q)$, and as $A$ can be taken the subgroup $(\GaL_5(q)/\{\pm1\})\times(\mF_q^* / (\mF_q^*)^2);$

\item $v=q^{16}$, $k=(q^3+1)(q^8-1)$, $\lambda=q^8+q^6-q^3-2$, $\mu=q^3(q^3+1)$,

\noindent $K\leq\AGaL_{16}(q)$, $\cG(K)$ is the affine half spin graph $\VD_{5,5}(q)$, and as $A$ can be taken the subgroup $(\mF_q^*\circ\operatorname{Inndiag}(D_5(q)):\aut(\mF_q)$.
\enmrt
Furthermore, every triple $(k,\lambda,\mu)$ of the attainable parameters for $K$ appears in items {\rm(i)--(ix)} and rows of Table~{\rm\ref{param2}} exactly once.
\ethrm

\begin{table}[ht]
\caption{Attainable parameters of small primitive affine rank 3 graphs}\label{param2}
\begin{tabular}{l|l|l|l|l|l|l}
\hline
$v=p^d$ & $k$ & $\lambda$ & $\mu$ & type & possible $A$ &  refs\\
\hline
$64=2^6$ & $18$ & $2$ & $6$ & C  & $3.\sym(6)$ & \cite[\S~10.24]{BrouwerM2022}\\
\hline
$169=13^2$ & $72$ & $31$ & $30$ & B  & $3\times (SL_2(3):4)$ & \cite[Table~11.8]{BrouwerM2022}\\
\hline
$243=3^5$ & $22$ & $1$ & $2$ & C  & $2\times M_{11}$ & \cite[\S~10.55]{BrouwerM2022}\\
\hline
$243=3^5$ & $110$ & $37$ & $60$ & C  & $2\times M_{11}$ & \cite[\S~10.55]{BrouwerM2022}\\
\hline
$256=2^8$ & $45$ & $16$ & $6$ & C  & $\sym(10)$ & \cite[\S~10.57]{BrouwerM2022}\\
\hline
$256=2^8$ & $102$ & $38$ & $42$ & C  & $L_2(17)$ & \cite[\S~10.58]{BrouwerM2022}\\
\hline
$361=19^2$ & $144$ & $59$ & $56$ & B  & $9\times\GL_2(3)$  & \cite[Table~11.8]{BrouwerM2022}\\
\hline
$625=5^4$ & $144$ & $43$ & $30$ & C  & $34.\sym(6)$ & \cite[\S~10.73A]{BrouwerM2022}\\
\hline
$625=5^4$ & $240$ & $95$ & $90$ & B  & $4.(2^4:\sym(6))$  & \cite[Table~11.8]{BrouwerM2022}\\
\hline
$841=29^2$ & $168$ & $47$ & $30$ & B  & $7\times(\SL_2(3):4)$  & \cite[Table~11.8]{BrouwerM2022}\\
\hline
$961=31^2$ & $240$ & $71$ & $56$ & B  & $15\times(2.\sym(4))$  & \cite[Table~11.8]{BrouwerM2022}\\
\hline
$961=31^2$ & $360$ & $139$ & $132$ & C  & $15\times\SL_2(5)$  & \cite[Table~11.8]{BrouwerM2022}\\
\hline
$1681=41^2$ & $480$ & $149$ & $132$ & C  & $40\circ\SL_2(5)$ & \cite[Table~11.7]{BrouwerM2022}\\
\hline
$2048=2^{11}$ & $276$ & $44$ & $36$ & C  & $M_{24}$ & \cite[\S~10.84]{BrouwerM2022}\\
\hline
$2048=2^{11}$ & $759$ & $310$ & $264$ & C  & $M_{24}$ & \cite[\S~10.85]{BrouwerM2022}\\
\hline
$2401=7^4$ & $240$ & $59$ & $20$ & C  & $6.S_4(3)$ & \cite[\S~10.89A]{BrouwerM2022}\\
\hline
$2401=7^4$ & $720$ & $229$ & $210$ & C  & $6.\sym(7)$ & \cite[\S~10.89C]{BrouwerM2022}\\
\hline
$2401=7^4$ & $960$ & $389$ & $380$ & C  & $48\circ\SL_2(5)$ & \cite[\S~10.89D]{BrouwerM2022}\\
\hline
$4096=2^{12}$ & $1575$ & $614$ & $600$ & C  & $(3\times HJ):2$ & \cite[\S~10.92]{BrouwerM2022}\\
\hline
$5041=71^2$ & $840$ & $179$ & $132$ & C  & $35\times\SL_2(5)$ & \cite[Table~11.7]{BrouwerM2022}\\
\hline
$6241=79^2$ & $1560$ & $419$ & $380$ & C  & $39\times\SL_2(5)$ & \cite[Table~11.7]{BrouwerM2022}\\
\hline
$6561=3^8$ & $1440$ & $351$ & $306$ & B  & $2.2^6:O_6^-(2).2$  & \cite[Table~11.8]{BrouwerM2022}\\
\hline
$15625=5^{6}$ & $7560$ & $3655$ & $3660$ & C  & $(2.HJ):4$ & \cite[\S~10.95]{BrouwerM2022}\\
\hline
$531441=3^{12}$ & $65520$ & $8559$ & $8010$ & C  & $2.\operatorname{Suz}.2$ & \cite[\S~10.100]{BrouwerM2022}\\
\hline
\end{tabular}
\end{table}

\prf Suppose first that $K\leq\sym(\Omega)$ is imprimitive, $\Sigma$ is a nontrivial system of blocks for $K$, and $\Delta\in\Sigma$ is one of the blocks. Recall that $K^\Delta$ denotes the permutation group on $\Delta$, induced by the action of the setwise stabilizer  of $\Delta$ in~$K$, and $K^\Sigma$ denotes the permutation group on $\Sigma$ induced by action of~$K$. Since $K$ is of rank~$3$, $\Sigma$ is the only nontrivial system of blocks, and the groups $K^\Delta$ and $K^\Sigma$ must be $2$-transitive. It follows that $\cG(K)$ is a disjoint union of cliques (the complement graph has a greater valency). Since $K$ has a normal regular subgroup, $K^\Delta$ also has such a subgroup. By Burnside's theorem, see, e.g., \cite[Theorem~4.1B]{DixM1996}, the socle of a $2$-transitive group $H$ is either a regular elementary abelian $p$-group for some prime $p$ or a nonregular nonabelian simple group. It follows that in the latter case, $H$ cannot have a normal regular subgroup, so $K^\Delta$ has a regular normal abelian elementary $p$-subgroup $U$, say, of order $p^t$. If the order of normal regular subgroup~$V$ of $K$ was divided by two distinct primes, then $\Sigma$ would not be the unique system of blocks for $K$, which is not the case. Thus, $|V|=p^{t+s}$, $K=V:A$ is affine, and $k$, $\lambda$, and $\mu$ are as in item~(i) of the theorem. A point stabilizer~$A$, considering as a subgroup of $GL(V)$, obviously stabilizes (setwise) the subspace~$U$. On the other hand, the setwise stabilizer of~$U$ in $GL(V)$ acts $2$-transitively on $U$ and the factor space $V/U$, so it can be taken as~$A$, which completes the proof of~(i).

Thus, we may assume that $K$ is primitive. The O'Nan--Scott theorem implies (see, e.g., \cite[Theorem~11.1.1]{BrouwerM2022}) that for every primitive rank $3$ group $K$, one of the following holds:
\nmrt[{\rm(a)}]
\item the socle of $K$ is a nonregular nonabelian simple group;
\item $K\leq K^\Delta\wr\sym(2)$ preserves a product decomposition $\Omega=\Delta^2;$
\item $K=V:A$ is affine, and $A\leq\GL(V)$ is a primitive linear group.
\enmrt

In case (a), $K$ has no any normal regular subgroup. In~(b), the group $K_0=K^\Delta$  must be $2$-transitive. As in the previous paragraph, $K^\Delta$ must include a regular normal elementary abelian $p$-subgroup $U$ for some prime $p$, say, of order $q=p^{t}$. Therefore, $V=U\times U$ is an elementary abelian $p$-group, proving (cf. (c)) that every group $K$ satisfying the hypothesis of the theorem is affine. Returning to (b), the group $\AGL_1(U)$ acts $2$-transitively on $U$, so we may take $\GL_1(U)$ as a point stabilizer in $K_0$. It remains to note that corresponding graph $\cG(K)$ is $q\times q$-grid, see \cite[Section~11.2]{BrouwerM2022}, so its parameters $k$, $\lambda$ and $\mu$ are as in (ii).

From now on, $K=V:A$ is a primitive affine group as in case (c) and $A$ is primitive linear group, since $K$ does not preserve a nontrivial product decomposition. The classification of the primitive affine rank 3 groups was completed by Liebeck in \cite{liebeckAffine} (see also \cite[Theorem~11.4.1]{BrouwerM2022}). He divided the affine groups of rank 3 into three classes: (A), (B), and~(C). Class (A) includes all the infinite series, while `extraspecial' class (B) and `exceptional' class (C) consist only of finitely many groups.

Suppose first that $K$ is not from (B) and (C), and $K\not\leq\AGaL_1(v)$. It follows from~\cite[Theorem~1.1]{SkrARS} that $K$ and $\cG(K)$ can be taken as in items (v)--(ix) of the theorem. The parameters $k$, $\lambda$ and $\mu$ can be easily extracted from \cite{BrouwerM2022}. It is worth noting two things here. First, the affine polar graph $\VO_{4}^{-}(q)$ and the Suzuki-Tits ovoid graph $\VSz(q)$, where $q=2^{2e+1}$, are not isomorphic, but have the same parameters $k$, $\lambda$, and $\mu$. Second, in item~(vii), we have to include the parameters of the complement of the affine polar graph $\VO_{2e}^{+}(2)$, because the valency of the later graph is greater than the valency of the complement.

If $K\leq\AGaL_1(v)$, then, up to a complement, $\cG(K)$ is either the Van Lint--Schrijver, Paley, or Peisert graph, see the short and clear proof of this fact in \cite[Theorem~3]{Muz21}, so the groups $K$ are described in~\cite[Propositions~1,2]{Muz21}, see items (iii)--(iv) of the theorem. The parameters of the corresponding graphs are taken from \cite{VLS} and \cite[Section~1.1.9]{BrouwerM2022}. The arguments from \cite{SkrARS} and \cite{GVW} allow us to avoid the coincidences with parameters of the graphs from items (v)--(ix) by setting $(p,e,t)\neq(2,3,2),(5,3,1),(2,3,3),(3,5,1),(2,5,2),(3,7,1),(2,11,1),(2,13,1)$ in the case of the Van Lint--Schrijver graphs. Note also that the parameters of the Paley and Peisert graphs of the same degree are always the same, though the graphs themselves are isomorphic only if $v=9$.

The affine rank $3$ groups $K$ from classes (B) and (C) are described in  \cite{liebeckAffine}, see also \cite[Tables~11.5 and 11.6]{BrouwerM2022}. The parameters of the corresponding graphs are extracted from \cite[Tables~11.7 and 11.8]{BrouwerM2022}. Let us emphasize that in Table~\ref{param2}, we list only the graphs from classes (B) and (C) that do not have the same parameters as the graphs listed in Theorem~\ref{180824c}.  One can verify this fact using the arguments of \cite[Section~3]{GVW} and Tables~4--6 therein.
\eprf

\begin{comment}
As an immediate consequence of Theorem~\ref{180824c}, we obtain the following statement to be used in Section~\ref{180824h}.

\crllrl{170824f}
The parameters of any rank~$3$ graph admitting a regular normal automorphism group are exactly those listed in Theorem~{\rm\ref{180824c}}.
\ecrllr
\end{comment}

\section{Proof of Theorem~\ref{100824a}}\label{180824h}

First, we  prove the sufficiency. The fact that the multivalued group $X(k)=\fX_{2k+1}(k)$ is coset if the number $4k+3$ is a prime power, was proved in~\cite{Ponomarenko2024a}. Let $\fX=\fX(\cG)$ be the multivalued group of a rank~$3$  graph $\cG$ with parameters $(k,\lambda,\mu)$  as  in Theorem~{\rm\ref{180824c}}.  Then the group $K=\aut(\cG)$ has a regular normal subgroup $G=V$. %by Corollary~\ref{170824f}.
By Theorem~\ref{170824q}, this implies that the multivalued group $\fX=\fX(\cG)=\fX(k,\lambda,\mu)$ is  isomorphic to a coset group $\fX(G,A,\id)$, where $A$ is a point stabilizer of~$K$.

To prove the necessity, let $\fX=(G,A,\phi) $ be  a multivalued  coset group for appropriate groups $G$, $A$, and a homomorphism~$\varphi$.  By the assumption, $\fX$ is defined on the set $\{e, x, y\}$ of the $A$-orbits.  We consider two cases depending on whether $x^\star$ is equal to $x$ or $y$.

Assume that $x^\star=y$. In this case, the Cayley digraph $\cay(G,x)$ is the Paley tournament, $|G|$ is a prime power, and the automorphism group of $\cay(G,x)$ has  a subgroup $A'$ acting semiregularly on  the nonidentity elements of~$G$ with  the same orbits as $A$ has, see \cite{Jones2020}. In particular,  we may assume that $A=A'$. Then
$$
n=|A|=|x|.
$$
On the other hand, from Theorem~\ref{170824a}, it follows  that $\fX=\fX_n(a)$, where $n=|A|$ and  $a=m_{x,x}^x$. Furthermore, by the assumption, the orbit $x$ contains no $g\in G$ such that $g^{-1}\in x$. This implies that
$$
n-2a=m_{x,y}^e={|A|\over{|x|}}=1.
$$
It follows that $\fX=X(k)$, where $k=a$. Moreover, since $2a=n-1$, we calculate $4k+3=4a+3=2n-2+3=2n+1=2|A|+1=2|x|+1=|G|.$ Thus,  $4k+3=|G|$ is a prime power, as required.

Now, let $x^\star=x$ and $y^\star=y$. It follows that the number $|x|$ (and also $|y|$) is even. Since $|x|$ divides $|A|$, this shows that the group $A$  is of even order. In addition,  the group $K=G\rtimes A$ is of  rank~$3$ and has the regular normal subgroup~$G$. It follows that the parameters $(k,\lambda,\mu)$ of the rank~$3$ graph $\cG=\cG(K)$ are  as in Theorem~\ref{180824c}. By Theorem~\ref{170824q},  we conclude that the multivalued groups~$\fX$ and $\fX(\cG)=\fX(k,\lambda,\mu)$ are isomorphic, and we are done.

\medskip

\textbf{Acknowledgements.} J. Guo and I. Ponomarenko were supported by the National Natural Science Foundation of China (Grant No. 12361003), A.~V.~Vasil'ev was supported by the Sobolev Institute of Mathematics state contract (FWNF-2022-0002) and by the National Natural Science Foundation of China (Grant No. 12171126).

%\bibliography{bibinp}{}
%\bibliographystyle{amsplain}
%\bibliography{/Users/inp1957/Yandex.Disk.localized/Inp/bibinp}{}

\end{document}